\documentclass{article}

\usepackage[top=1.5in,left=1.25in,right=1.25in,bottom=1.5in]{geometry}

\usepackage{latexsym}
\usepackage[pdftex]{graphicx}
\usepackage{epsfig, ecltree, epic, eepic}
\usepackage{enumerate,amsmath,amssymb,dsfont,pifont,amsthm}
\usepackage{xcolor}
\usepackage{ifthen}
\usepackage{graphicx}
\usepackage[arrow, matrix, curve]{xy}
\usepackage{amsthm,amsmath}
\usepackage{stmaryrd}
\usepackage{comment}
\usepackage{enumerate,amssymb,dsfont,pifont}
\usepackage{tikz}
\usepackage{hyperref}
\hypersetup{
    colorlinks,
    citecolor=black,
    filecolor=black,
    linkcolor=black,
    urlcolor=black
}

\theoremstyle{plain}
\newtheorem{theorem}{Theorem}[section]

\newtheorem{lemma}[theorem]{Lemma}
\newtheorem{proposition}[theorem]{Proposition}

\theoremstyle{definition}
\newtheorem{definition}[theorem]{Definition}

\theoremstyle{remark}
\newtheorem{remark}[theorem]{Remark}

\newtheorem{example}[theorem]{Example}

\newcommand{\bbc}{\mathbb{C}}
\newcommand{\bbr}{\mathbb{R}}

\newcommand{\bbp}{\mathbb{P}}
\newcommand{\bbe}{\mathbb{E}}

\newcommand{\bbn}{\mathbb{N}}

\newcommand{\cf}{\mathcal{F}}
\newcommand{\cm}{\mathcal{M}}

\newcommand{\hx}{\hat{X}}
\newcommand{\ip}{\bbp^{\tau,x}}

\newcommand{\abs}[1]{\left| #1 \right|}
\newcommand{\norm}[1]{\left\| #1 \right\|}

\newcommand{\roi}[2]{\left[ #1 , #2 \right[}

\newcommand{\ci}[2]{\left[ #1, #2 \right]}

\newcommand{\rosi}[2]{\roi{[#1}{#2[}}

\newcommand{\csi}[2]{\ci{[#1}{#2]}}

\makeatletter
\newcommand*\Si{\mathpalette\Si@{.5}}
\newcommand*\Si@[2]{\mathbin{\vcenter{\hbox{\scalebox{#2}{$\m@th#1\bullet$}}}}}
\makeatother

\newcommand{\ls}{\llbracket}
\newcommand{\rs}{\rrbracket}

\renewcommand{\leq}{\leqslant}
\renewcommand{\geq}{\geqslant}

\begin{document}

\newcounter{farbe}
    \setcounter{farbe}{1} 

\ifthenelse{\value{farbe}=1}{\newcommand{\cor}{\color{blue}}}{\newcommand{\cor}{\color{black}}}

\allowdisplaybreaks

\title{\bfseries The Time-Dependent Symbol of a Non-Homogeneous Itô Process and corresponding Maximal Inequalities}

\author{%
    \textsc{Sebastian Rickelhoff}%
    \thanks{University of Siegen, Department of Mathematics, Walter-Flex-Straße 3,
              D-57072 Siegen, Germany,
              \texttt{Sebastian.Rickelhoff@uni-siegen.de}, 
              phone: +49 (271) 740-3575.} \ and
                  \textsc{Alexander Schnurr}%
    \thanks{University of Siegen, Department of Mathematics,Walter-Flex-Straße 3,
               D-57072 Siegen, Germany,
              \texttt{Schnurr@mathematik.uni-siegen.de}, 
              phone: +49 (271) 740-3806.}
    }

\date{\today}

\maketitle
\begin{abstract} 
\noindent
The probabilistic symbol is defined as the right-hand side derivative at time zero of the characteristic functions corresponding to the one-dimensional marginals of a time-homogeneous stochastic process. As described in various contributions to this topic, the symbol contains crucial information concerning the process. 
When leaving time-homogeneity behind, a modification of the symbol by inserting a time component is needed.\\
In the present article we show the existence of such a time-dependent symbol for non-homogeneous Itô processes. Moreover, for this class of processes we derive maximal inequalities which we apply to generalize the Blumenthal-Getoor indices to the non-homogeneous case. These are utilized to derive several properties regarding the paths of the process, including the asymptotic behavior of the sample patsh, the existence of exponential moments and the finiteness of p-variation.  In contrast to many situations where non-homogeneous Markov processes are involved, the space-time process \emph{cannot} be utilized when considering maximal inequalities.
\end{abstract}



\section{Introduction}
The \textit{probabilistic symbol} $p(x,\xi)$ of a Markov process $X$ is the function $p:\bbr^d\times \bbr^d\to \bbc$ given by
\begin{align} \label{stoppedsymbol} 
     p(x,\xi):=- \lim_{t\downarrow 0}\bbe^x \frac{e^{i(X^\sigma_t-x)'\xi}-1}{t}
\end{align}
if the limit exists and coincides for every $R>0$, where
$$
\sigma:=\inf\{t\geq 0 : \lVert X_t^x-x \rVert > R \}.
$$
This symbol proves to be a crucial concept for deriving a wide range of properties of the stochastic process,  such as conservativeness (cf. \cite{schilling98pos}, Theorem 5.5),  asymptotic behavior (cf.  \cite{generalizedindices},Theorems 3.11 and 3.12), strong $\gamma$-variation (cf. \cite{sdesymbol} Corollary 5.10), Hausdorff-dimension (cf. \cite{schilling98hdd}, Theorem 4) and H\"older conditions \cite{schilling98}. For a survey on recent results we refer to \cite{levymatters3} and \cite{hope}.  By now, all of these results are restricted to the time-homogeneous case. \\
Proving some of these results, the symbol is utilized to derive \textit{maximal-inequalities}.  Inequalities of this kind have been proved for L\'evy processes (cf. \cite{pruitt}), certain Feller processes (cf. \cite{schilling98}) and homogeneous diffusions with jumps (cf. \cite{generalizedindices}). One can find a throughout discussion of maximal inequalities for various classes of stochastic processes in \cite{MaxInSchKueh}.
Formulating these maximal inequalities,  \textit{Blumenthal-Getoor indices} are used (cf.  \cite{schilling98},\cite{generalizedindices}) which allow for a governance of the process's paths by the behavior of the symbol in the variable $\xi$, and, therefore, for the derivation of the properties stated above.

However, when leaving the time-homogeneity behind one does not expect the symbol, being the right-hand side derivative of the characteristic functions corresponding to the one-dimensional marginal \textit{at time zero},  to yield any information regarding the entire process. To overcome this,  \cite{rueschendorf} proposed adding a time component to the symbol or more precisely
\begin{equation}\label{intro:timedep-sym}
p(\tau,x,\xi) := -\lim_{h \downarrow 0} \bbe^{\tau,x} \left( \frac{e^{i(X^{\sigma}_{\tau+h}-x)'\xi}-1}{h} \right) \text{ for } x,\xi\in\bbr^d, \ \tau\geq 0.
\end{equation}
Moreover,  the existence of such a \textit{time-dependent probabilistic symbol} was shown for rich càdlàg Feller evolution processes, i.e.,   non-homogeneous càdlàg Markov processes such that
$$
T_{\tau,t}u(x):= \bbe^{\tau,x}u(X_t)
$$
for $0\leq \tau \leq t$ and $u \in B_b(\bbr^d)$ forms a strongly continuous evolution system. In addition, the domain of the \textit{infinitesimal generator} $A_\tau$
$$
A_\tau f:= \lim_{h \downarrow 0} \frac{T_{\tau,\tau+h}f-f}{h}, \quad \tau>0
$$
given by 
$$
D(A_\tau):= \left\{ f  : \lim_{h \downarrow 0} \frac{T_{\tau,\tau+h}f-f}{h} \text{ exists } \right\}
$$
contains the test functions $C_c(\bbr^d)$. Theorem 4.5 of \cite{rueschendorf} shows that the generator $A_\tau|_{C_c^\infty}$ is a \textit{pseudo-differential operator} with symbol $-q(\tau,x,\xi)$, i.e., 
$$
A_\tau f(x) = - \frac{1}{(2\pi)^{\frac{d}{2}}} \int_{\bbr^d} e^{ix'\xi}q(\tau,x,\xi)\hat{f}(\xi) \ d\xi
$$
for $\tau \geq 0, x,\xi\in \bbr^d$ and $\hat{f}$ is the Fourier transform of $f$.  Moreover, it is shown that the symbol of the generator $q(\tau,x,\xi)$ and the time-dependent probabilistic symbol defined in \eqref{intro:timedep-sym} $p(\tau,x,\xi)$ coincide if the symbol is continuous in $x$. \\
In the present article,  we prove the existence of the time-dependent symbol for non-homogeneous Itô processes (cf. Definition \ref{non-hom ito}).  We utilize this result to prove maximal-inequalities,  the existence of non-homogeneous generalizations of the Blumenthal-Getoor indices and an exemplary selection of properties of such processes.  Before we do so, we fix some notations:\\
A family of $\sigma$-fields $(\mathcal{G}_t^\tau)_{0 \leq \tau \leq t}$ is called a \textit{two-parameter filtration} if $\mathcal{G}_s^\tau \subset \mathcal{G}_t^\tau$ for all $0 \leq \tau \leq s \leq t$ and $\mathcal{G}_t^{\tau_1} \subset \mathcal{G}_t^{\tau_2}$ for $0\leq \tau_1 \leq \tau_2 \leq t$. 
The \textit{natural double filtration} of $X$ is denoted by $(\mathcal{F}^X)_t^\tau)_{0\leq \tau \leq t}$ and is defined as
$$
\left(\mathcal{F}^X\right)_t^\tau:= \sigma(X_s: \tau \leq s \leq t).
$$
Let $(\Omega,\mathcal{M})$ be a measurable space equipped with the two-parameter filtration $ (\mathcal{M}_t^\tau)_{0\leq \tau \leq t}$. We call a stochastic process $X$ \textit{adapted to the two-parameter filtration} if for all $0 \leq \tau \leq t$ 
$$
\left(\mathcal{F}^X\right)_t^\tau \subset \mathcal{M}_t^\tau.
$$
We tacitly assume that every stochastic process $X:=(X_t)_{t \geq 0}$ is defined on a generic stochastic basis $(\Omega, \mathcal{A}, (\mathcal{A}_t)_{t \geq 0},\bbp)$ takes values in $(\bbr^d, \mathcal{B}(\bbr^d))$ and is c\'{a}dl\'{a}g. Here, $\mathcal{B}(\bbr^d)$ is the $\sigma$-field of Lebesgue sets.  Moreover, we call $\Delta X_t := X_t - \lim_{s \uparrow t} X_s$ the jump of the process at time $t\geq 0$, and for a stopping time $\tau$ we call $X^\tau:= X1_{\ls 0, \tau \rs}+ X_\tau 1_{\ls \tau, \infty \ls}$ the stopped process. The stochastic interval $\ls \tau , \sigma \ls$ for two stopping times $\tau, \sigma$ is defined by $\{ (\omega, t) \in \Omega\times \bbr_+: \tau(\omega) \leq t < \sigma(\omega)\}$. The stochastic intervals $\ls \tau , \sigma \rs$, $\rs \tau , \sigma \ls$, $\rs \tau , \sigma \rs$ are defined alike. \\
A (strong) \emph{Markov process} $(\Omega, \mathcal{M}, ( \mathcal{M}_t^\tau)_{0\leq \tau \leq t}, (X_t)_{t \geq 0}, \bbp^{\tau,x})_{\tau\in \bbr_+, x \in \bbr^d}$ satisfies 
\begin{equation}
\bbe^{\tau,x} \left[ f(X_t) \mid \mathcal{M}_s^\tau \right] = \bbe^{s,X_s}\left( f(X_t \right)) ,\quad \bbp^{\tau,x}\text{-a.s.}
\end{equation}
for all $\tau \leq s \leq t$ and all bounded Borel-measurable functions $f$.  Moreover, every Markov process is \textit{normal}, i.e.,  $\bbp^{\tau,x}(X_\tau=x)=1$. For a more information on Markov processes see \cite{gulcast} and \cite{cast}.
We associate an evolution system $(T_{\tau,t})_{0\leq \tau \leq t}$ of operators on $B_b(\bbr^d)$ with every Markov process by setting
\[
    T_{\tau,t} u(x):= \bbe^{\tau,x} u(X_t).
\]

\section{The Time-Dependent Probabilistic Symbol}
In this section $X:=(\Omega, \mathcal{M}, ( \mathcal{M}_t^\tau)_{0\leq \tau \leq t}, (X_t)_{t \geq 0}, \bbp^{\tau,x})_{\tau\in \bbr_+, x \in \bbr^d}$ denotes a Markov process. 
Before we start with the main topic of this section, we properly define the space-time process of a Markov process. That is due to the fact that various different definitions are used in the literature.  The following definitions follows \cite{boettcher}. \\
Let $\hat{\Omega}:= \bbr_+ \times \Omega$ and the $\sigma$-field $\hat{\cm}:= \{ A \subset \hat{\Omega}: A_s \in \cm \ \forall
s \in \bbr_+\}$ where $A_s$ denotes the $s$-slice of $A$ for $s \geq 0$.
We define a process $\hat{X}$ with values in $\bbr_+ \times \bbr^d$ by
\begin{equation*}
\hat{X}_t(\hat{\omega})= \hat{X}_t((c, \omega)):= (t+c, X_{t+c}(\omega)).
\end{equation*}
Moreover, we set
\begin{equation*}
\hat{\theta}_\tau: \hat{\Omega} \to \hat{\Omega}; (s,\omega) \mapsto (s+\tau, \omega)
\end{equation*}
 and 
\begin{equation}\label{space-time prob}
\bbp^{(\tau,x)}(A):= \bbp^{\tau,x} \left( \pi_0^{-1} \left( \hat{\theta}_\tau^{-1} (A) \right) \right)
\end{equation}
where $\pi_0 :\Omega \to \hat{\Omega}; \omega \mapsto (0, \omega)$, and $A \in \hat{\cm}$. \\
We call the homogeneous Markov process 
$$
\hat{X}:=(\hat{\Omega}, \hat{\cf}, (\cf_t^\tau)_{t\geq \tau}, (\hat{X}_t)_{t\geq 0}, (\hat{\theta}_t)_{t \geq 0},\hat{\bbp}^{(\tau,x)})_{(\tau,x)\in \bbr_+ \times \bbr^d}
$$
the \textit{space-time process associated with $X$}.  The transition probability function of $X$ is given by 
\begin{equation*}
\hat{P}(t,(\tau,x), A):= P(\tau,x; t+\tau, A_{t+ \tau}),
\end{equation*}
and for any $\hat{\mathcal{F}}$-measurable random variable $Y$ it holds true that
$$
Y=Y \circ \hat{\theta}_\tau \circ \pi_0, \quad \bbp^{(\tau,x)}\text{-a.s.},
$$
i.e.,  for $\bbp^{(\tau,x)}$-almost all $(c,\omega) \in \hat{\Omega}$ we have 
\begin{equation}
Y(c,\omega)= Y(\tau, \omega).
\end{equation}
This will be used frequently throughout the following calculations.
\begin{definition}\label{inhomMarkSem}
We call a Markov process $X$ \textit{non-homogeneuous Markov semimartingale} if for every $\bbp^{\tau,x}$, $\tau \geq 0, x \in \bbr^d$ the process $(X_t)_{t \geq 0}$ is a semimartingale on $[\tau, \infty)$.
\end{definition}

\begin{definition}\label{non-hom ito}
We call a non-homogeneous Markov semimartingale \textit{non-homogeneous Itô process} if its characteristics $(B,C,\nu)$ are of the form
\begin{align}\label{non-hom dif jump}
B_t^{(i)} &= \int_\tau^t b^{(i)}(s,X_s) \ ds \quad \bbp^{\tau,x}\text{-a.s.}  \nonumber\\
C_t^{(ij)} &= \int_\tau^t c^{(ij)}(s,X_s) \ ds \quad \bbp^{\tau,x}\text{-a.s.}\\
\nu( ; dt,dx) &= dt N_t(X_t,dx) \quad \bbp^{\tau,x}\text{-a.s.} \nonumber
\end{align}
for $t \geq \tau \geq 0, x \in E$ and $i,j \in \{1,...,d\}$.
\end{definition}

\begin{definition}
Let $X$ be a Markov process and let 
\begin{equation}\label{exittimeR}
\sigma:= \sigma_R^{\tau,x}:= \inf\{ h \geq \tau: \Vert X^{\tau,x}_{h} -x \Vert > R\}=\tau+\inf\{ h \geq 0: \Vert X^{\tau,x}_{\tau+h} -x \Vert > R\}
\end{equation}
be the first exit time from the ball of radius $R>0$ after $\tau\geq 0$, and $\Vert \cdot \Vert$ the maximum norm. The function $p:\bbr_+\times \bbr^d \times \bbr^d \to \bbc$ defined by 
\begin{equation}\label{inhom-symbol}
p(\tau,x,\xi) := -\lim_{h \downarrow 0} \bbe^{\tau,x} \left( \frac{e^{i(X^\sigma_{\tau+h}-x)'\xi}-1}{h} \right)
\end{equation}
is called the \textit{time-dependent probabilistic symbol} of the process, if the limit exists for every $\tau \geq 0$ and $x,\xi \in \bbr^d$ independently of the choice of $R$.
\end{definition}

\begin{example}\label{ex: time-dep symbol}
\begin{itemize}
\item[(a.)]
Let $(X_t)_{t \geq 0}$ be an additive process on $(\Omega, \mathcal{F}, (\mathcal{F}_t)_{t\geq 0},\bbp)$ in the sense of Definition 1.6 of \cite{sato}, i.e.,  $X$ has independent increments, is stochastically continuous and càdlàg and starts in $0$ at time $t=0$.  We define a family of probability measures on $(\Omega, \mathcal{F})$ by 
$$
\bbp^{\tau,x}(X_t \in B) := \bbp( X_t -X_\tau \in B-x) 
$$
for $B \in \mathcal{E}$ and $\tau \in \bbr_+, x \in E$.  For this family it holds true that $\bbp^{\tau,x}(X_\tau=x)=1$, and $(\Omega, \mathcal{F}, (\mathcal{F}_t)_{t\geq 0},(X_t)_{t \geq 0} ,\bbp^{\tau,x})_{\tau \geq 0, x \in E}$ is a Markov process.  \\
Let, in addition, $X$ be a semimartingale for all $\bbp^{\tau,x}$ which is quasi-left-continuous, and possesses the characteristics $(B,C,\nu)$. Theorem II.4.15 of \cite{jacodshir} provides the existence of a version of $(B,C,\nu)$ that is deterministic. Hence, in the following we assume $(B,C,\nu)$ to be deterministic. By Corollary II.4.18 of \cite{jacodshir} $X$ has no fixes times of discontinuity, i.e.,  
$$\{t\geq 0: \nu(\{t\} \times E)>0\} = \emptyset.$$
Therefore, when calculating the time-dependent probabilistic symbol of $X$ we derive with Theorem II.4.15 of \cite{jacodshir} that
\begin{align*}
p(\tau,x,\xi) = &-\lim_{h \downarrow 0} \bbe^{\tau,x} \left( \frac{e^{i(X_{\tau+h}-x)'\xi}-1}{h} \right) \\
= &-\lim_{h \downarrow 0} \bbe \left( \frac{e^{i(X_{\tau+h}-X_\tau)'\xi}-1}{h} \right) \\
=& -\lim_{h \downarrow 0} \frac{e^{ i \xi (B_{\tau+h} -B_\tau) - \frac{1}{2}\xi' ( C_{\tau+h}-C_\tau)\xi + \int_{\bbr^d \setminus \{0\}} e^{i\xi'y}-1-i\xi'y\chi(y) \ \nu ((\tau,\tau+h],dy) }-1}{h}.
\end{align*}
This limit exists if and only if 
$$
-\lim_{h \downarrow 0} \frac{ i \xi (B_{\tau+h} -B_\tau) - \frac{1}{2}\xi' ( C_{\tau+h}-C_\tau)\xi + \int_{\bbr^d \setminus \{0\}} e^{i\xi'y}-1-i\xi'y\chi(y) \ \nu ((\tau,\tau+h],dy)}{h}
$$
exists.\\
If $B^{(i)}$ and $C^{(ij)}$ are right-differentiable for all $i,j \in \{1,...,d\}$ and if the function 
$$
\tau \mapsto \int_{\bbr^d \setminus \{0\}} e^{i\xi'y}-1-i\xi'y\chi(y) \ \nu ((0,\tau],dy)
$$
is right-differentiable the time-dependent symbol exists and is of the form 
\begin{equation*}
p(\tau,x,\xi)=  i \xi \partial_+ B_{\tau} - \frac{1}{2}\xi' \partial_+ C_{\tau} \xi + \partial_+ \int_{0}^\tau \int_{\bbr^d \setminus \{0\}} e^{i\xi'y}-1-i\xi'y\chi(y) \ \nu (ds,dy).
\end{equation*}

\item[(b.)] Let $(X_t)_{t \geq 0}$ be a one-dimensional Brownian motion with variance function $\sigma^2(t)$ on \\$(\Omega, \mathcal{F}, (\mathcal{F}_t)_{t\geq 0},\bbp)$.  The process $X$ is additive and a continuous semimartingale and 
$$
\bbe\left( e^{i(X_{t+\tau}-X_\tau)'\xi} \right)= \exp\left( - \frac{1}{2}(\sigma^2(t+\tau)-\sigma^2(\tau))\xi^2 \right).
$$
By the previous example the (non-homogeneous) probabilistic symbol exists if and only if the variance function is right-differentiable with right-derivative $\partial_+ \sigma^2$.  In this case we have 
$$
p(\tau,x,\xi) = -\frac{1}{2}\xi^2 \partial_+ \sigma^2(\tau).
$$
\end{itemize}
\end{example}

We have seen in \eqref{intro:timedep-sym} that under some mild conditions for a rich càdlàg Feller evolution process the symbol of the generator $q(\tau,x,\xi)$ and the time-dependent probabilistic symbol $p(\tau,x,\xi)$ coincide. Additionally, Corollary 3.5 of \cite{boettcher} states that the symbol of the generator of the homogeneous space-time process $\hx$ corresponding to $X$ is given by 
$$
q(\hat{x},\hat{\xi}) = -i\xi_0 + q(\tau,x,\xi)
$$
with $\tau \geq 0, x \in \bbr^d, \xi \in \bbr^d$, $\hat{x}=(\tau,x),\hat{\xi}=(\xi_0,\xi) \in \bbr^{d+1}$. Therefore,  we expect the space-time process to be useful when calculating the symbol of a non-homogeneous Itô process,  provided the characteristics of the space-time process are the ones of a homogeneous diffusion with jumps.\\
Hence, the following lemma states the characteristics of the space-time process associated to a non-homogeneous Markov semimartingale. Note that the proof is omitted since it is quite straightforward but needs tedious calculations.

\begin{lemma}\label{lem:space-time sem}
Let $X$ be a non-homogeneous Markov semimartingale with characteristics $(B,C,\nu)$.  In this case the space-time process $\hat{X}$ associated with $X$ is semimartingale for all $\bbp^{(\tau,x)}, (\tau,x)\in \bbr_+ \times \bbr^d$ and its characteristics $(\hat{B},\hat{C},\hat{\nu})$ are given by
\begin{align}
\hat{B}_t(c, \omega)&= (t,B_{\tau+t}(\omega)), \quad \bbp^{(\tau,x)}\text{-a.s.}\label{CharSpace-Timea}\\
\hat{C}_t(c, \omega)&= \begin{pmatrix}
  0 & \begin{matrix} \cdots & \hspace{-1mm} 0 \end{matrix} \\[-0.5em]
  \begin{matrix} \vdots \\ 0 \end{matrix}  &
  \begin{matrix}\\
 \raisebox{\height}{$C_{\tau+t}(\omega)$} 
  \end{matrix}
\end{pmatrix},  \quad  \bbp^{(\tau,x)}\text{-a.s.}\label{CharSpace-Timeb}\\
\hat{\nu}((c,\omega); ds,  du \times dy) &= \nu(\omega; ds+\tau, dy) \delta_0(du),  \quad  \bbp^{(\tau,x)}\text{-a.s.}\label{CharSpace-Timec}
\end{align}
for $\tau \leq t$ and $(c,\omega) \in \hat{\Omega}$.
\end{lemma}

Before we state the main theorem of this section, the following lemma provides that non-homogeneous Itô process are indeed a generalization of rich Feller evolution processes:

\begin{lemma}
Let $(X_t)_{t\geq 0}$ be a rich Feller evolution process on $C_\infty(\bbr^d)$ with generator $A_s$ and time dependent symbol $p: \bbr_+ \times \bbr^d \times \bbr^d \to \bbc$ such that $p(\cdot, x,\xi)$ is continuous for all $x, \xi \in \bbr^d$. Then $X$ is a non-homogeneous Itô process.
\end{lemma}
\begin{proof}
Let $(X_t)_{t\geq 0}$ be a rich Feller evolution system on $C_\infty(\bbr^d)$. Analogously to the proof of Theorem 3.1 of \cite{mydiss}, but under usage of the non-homogeneous version of Dynkin's formula as mentioned in \cite{rueschendorf}, we show that $(X_t)_{t\geq 0}$ is a semimartingale with characteristics $(B,C,\nu)$. Since it is well-known that $X$ is a Markov process, it suffices to show that the characteristics are of the form mentioned in Definition \ref{non-hom ito}: \\
By Theorem 3.2 and Lemma 3.7 of \cite{boettcher} the space-time process $\hat{X}$ associated to $X$ is a rich Feller process on $C_\infty(\bbr_+ \times \bbr^d)$.  Hence,  Theorem 3.10 of \cite{mydiss} provides that the characteristics $(\hat{B},\hat{C},\hat{\nu})$ of $\hat{X}$ are of the form 
\begin{align*}
\hat{B}_t^{(i)} &= \int_0^t b^{(i)}(\hat{X}_s) \ ds, \quad \bbp^{(\tau,x)}\text{-a.s.}  \nonumber\\
\hat{C}_t^{(ij)} &= \int_0^t c^{(ij)}(\hat{X}_s) \ ds, \quad \bbp^{(\tau,x)}\text{-a.s.}\\
\hat{\nu}(\;  ; dt,d\hat{x}) &= dt N(\hat{X}_t,d\hat{x}), \quad \bbp^{(\tau,x)}\text{-a.s.} \nonumber
\end{align*}
and, therefore,  we conclude with equations \eqref{CharSpace-Timea} to \eqref{CharSpace-Timec} that 
\begin{align*}
B_t^{(i)} &= \int_0^{t-\tau} b^{(i)}(\hat{X}_s) \ ds= \int_0^{t-\tau} b^{(i)}(s+\tau,X_{s+\tau}) \ ds= \int_\tau^{t} b^{(i)}(s,X_{s}) \ ds, \quad \bbp^{(\tau,x)}\text{-a.s.}  \\
C_t^{(ij)} &= \int_0^{t-\tau} c^{(ij)}(\hat{X}_s) \ ds= \int_0^{t-\tau} c^{(ij)}(s+\tau,X_{s+\tau}) \ ds = \int_\tau^{t} c^{(ij)}(s,X_{s}) \ ds, \quad \bbp^{(\tau,x)}\text{-a.s.}\\
\nu( \; ; dt,dx) &= dt N(X_{t},dx) \quad \bbp^{(\tau,x)}\text{-a.s.} 
\end{align*}
for $t \geq \tau$.
Transition from $ \bbp^{(\tau,x)}$ to $ \bbp^{\tau,x}$ yields the statement. 
\end{proof}

\begin{theorem}\label{main:symbol}
Let $X$ be a non-homogeneous Itô process and let $\ell= (\ell^{(j)})_{1 \leq j \leq d}$ and $Q= (q^{(ik)})_{1 \leq j,k\leq d}$ be continuous and locally bounded, $N$ be such that the function 
$$
(s,x) \mapsto \int_{y \neq 0} (1 \wedge y^2) \ N_s(x,dy)
$$ 
is continuous and locally bounded. In this case the time-dependent symbol exists and equals
$$
p(\tau,x,\xi) = -i \ell(\tau,x)'\xi + \frac{1}{2} \xi' Q(\tau,x) \xi - \int_{y \neq 0} \left(e^{iy' \xi} -1 - iy' \xi \cdot \chi(y) \right) \ N_\tau(x,dy).
$$
\end{theorem}
\begin{proof}
Let $X$ be a non-homogeneous Itô process, and let $\hx$ be the associated space-time process.  The characteristics of  $\hat{X}$ are given by 
\begin{align*}
\hat{B}_t^{(1)}&= t = \int_0^t 1 \ ds\\
\hat{B}_t^{(i)} &= \int_\tau^{t+\tau} \ell^{(i-1)}(s,X_s) \ ds =\int_0^{t} \ell^{(i-1)}(s+\tau,X_{s+\tau}) \ ds =\int_0^{t} \ell^{(i-1)}(\hx_{s}) \ ds ,\\
\hat{C}_t^{(1j)} &= \int_0^{t} 0 \ ds \\
\hat{C}_t^{(ij)} &= \int_0^{t} c^{((i-1,j-1)}(\hx_{s}) \ ds 
\end{align*}
for $i,j \in \{2,...,d+1\}$. All of the equations above are meant $\hat{\bbp}^{(\tau,x)} \text{-a.s.}$. For $T \in \mathcal{B}(\bbr_+)$ and $B \in \mathcal{B}(\hat{E})$ we have
\begin{align*}
\hat{\nu}( \; ; T,B) &= \nu(\;  ; T+\tau, B_0) =  \int_{T+\tau} N_{s}(X_{s},B_0)\ ds\\
&=\int_{T} N_{s+\tau}(X_{s+\tau},B_0)\ ds =\int_{T} \hat{N}((s+\tau,X_{s+\tau}),B)\ ds \\
&=\int_{T} \hat{N}(\hx_s,B)\ ds, \quad \bbp^{(\tau,x)}\text{-a.s.}
\end{align*}
where $\hat{N}((c,\omega),B)=N_c(\omega,B_0)$ is a transition kernel from $\hat{\Omega}\times \hat{E}$ to $\bbr^d$.  Consequently,  $\hx$ is a Itô process and we denote by $\hat{p}(\hat{x},\hat{\xi})$ its (homogeneous) probabilistic symbol.  \\
For the stopping time $\hat{\sigma}:= \inf\{ h \geq 0: \Vert \hat{X}_h -(\tau,x) \Vert > R \}$ it holds true that 
\begin{align*}
\hat{\sigma}(\hat{\theta}_\tau(\pi_0(\omega)))&:= \inf\{ h \geq 0: \Vert \hat{X}_h(\hat{\theta}_\tau(\pi_0(\omega))) -(\tau,x) \Vert > R \}\\
&= \inf\{ h \geq 0: \Vert \hat{X}_h(\tau,\omega) -(\tau,x) \Vert > R \}\\
&= \inf\{ h \geq 0: \Vert (\tau+h,X_{\tau+h}(\omega)) -(\tau,x) \Vert > R \}\\
&= \inf\{ h \geq 0: \Vert (h,X_{\tau+h}(\omega)-x) \Vert > R \}\\
&= R \wedge \inf\{ h \geq 0: \Vert X_{\tau+h}(\omega)-x \Vert > R \}\\
&= R \wedge ( \sigma -\tau).
\end{align*}
We compute for $\hat{\xi}=(0,\xi)$, $\hat{\ell}=(1,\ell), \hat{c}=(\hat{c}^{(ij)}$ with $\hat{c}^{(1,j)}=\hat{c}^{(i,1)}=0$ and $\hat{c}^{(ij)}=c^{((i-1),(j-1))}$ for $i,j \in \{1,...,d+1\}$.
\begin{align*}
p(\tau,x,\xi) 
&=  -\lim_{h \downarrow 0} \bbe^{\tau,x} \left(\frac{e^{i (X^\sigma_{\tau+h} -x)'\xi}-1}{h} \right)\\
&=  -\lim_{h \downarrow 0} \bbe^{\tau,x} \left(\frac{e^{i (X_{\tau+(h\wedge (\sigma-\tau))} -x)'\xi}-1}{h} \right)\\
&=  -\lim_{\underset{h <R}{h \downarrow 0}} \bbe^{\tau,x} \left(\frac{e^{i (X_{\tau+(h\wedge (\sigma-\tau)\wedge R)} -x)'\xi}-1}{h} \right)\\
&=  -\lim_{\underset{h <R}{h \downarrow 0}} \bbe^{\tau,x} \left(\frac{e^{i (\hx_{h\wedge (\hat{\sigma}\circ \hat{\theta}_\tau \circ \pi_0)} \circ \hat{\theta}_\tau \circ \pi_0 -(\tau,x))'\hat{\xi}}-1}{h} \right)\\
&= -\lim_{\underset{h <R}{h \downarrow 0}} \hat{\bbe}^{(\tau,x)} \left(\frac{e^{i (\hx^{\hat{\sigma}}_h -(\tau,x))'\hat{\xi}}-1}{h} \right)\\
&= \hat{p}((\tau,x), \hat{\xi})\\
&= -i \ell(\tau,x)'\hat{\xi}+ \frac{1}{2} \hat{\xi}' c(\tau,x) \hat{\xi}- \int_{\hat{y}\neq 0} \left(e^{i\hat{y}'\hat{\xi}}-1-i \hat{y}' \hat{\xi}\chi(\hat{y}) \right) \ \hat{N}((\tau,x),dy)\delta_0(du)\\
&= -i \ell(\tau,x)'\xi+ \frac{1}{2} \xi' c(\tau,x) \xi- \int_{y\neq 0} \left(e^{iy'\xi}-1-i y' \xi\chi(y) \right) \ N_\tau(x,dy).
\end{align*}
\end{proof}

We have observed, for instance in Example \ref{ex: time-dep symbol} that a process does not necessary need to be a non-homogeneous Itô process to possess a time-dependent symbol.  However, the processes considered up to this point are quasi-left continuous.  We will see in the subsequent example that, when leaving quasi-left continuity behind, the time-dependent symbol does not contain the same information on the process as before.  This is not unexpected, as we have encountered similar situations in the homogeneous case.

\begin{example} \label{example:det3}
Let us consider the following (deterministic!) example:
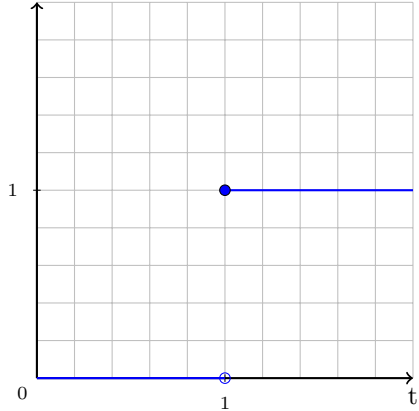
\begin{center}
\begin{tikzpicture}[x=.5cm, y=.5cm,domain=-9:9,smooth]
   \draw [color=gray!50]  [step=5mm] (0,0) grid (10,10);
   \draw[->,thick] (0,0) -- (10,0) node[below]{t};
   \draw[->,thick] (0,0) -- (0,10) node[above]{} ;
  
      \draw (5,-.1) -- (5,.1) node[below=5pt] {$\scriptstyle 1$};

     \draw (-.1,5) -- (.1,5) node[left=5pt] {$\scriptstyle 1$};

   \node[below left]{$\scriptstyle0$};
\draw[color=blue,thick] (0,0) -- (5,0);
\draw[color=blue,thick] (5,5) -- (10,5);
\draw[fill=blue] (5,5) circle(0.7mm) ;
\draw[color=blue] (5,0) circle(0.7mm) ;
\end{tikzpicture} 
\end{center}
By adding all other starting points in a Markovian manner (c.f.  \cite{detmp2}), and considering the truncation function $h(x)=x1_{\{|x| \leq 0.5\}}$ we receive a non-homogeneous Markov semimartingale with characteristics $B \equiv 0, C \equiv 0$ and $\nu(dt,dx)= \delta_1(dt)\delta_1(dx)$.  In this case, the time-dependent symbol exists and is given by 
$$
p(\tau,x,\xi) =  -\lim_{h \downarrow 0} \bbe^{\tau,x} \left(\frac{e^{i (X^\sigma_{\tau+h} -x)'\xi}-1}{h} \right) \equiv 0.
$$
Nevertheless, the symbol does not provide any information regarding the process.
\end{example}

\section{Maximal Inequalitys and Time-Dependent Blumenthal-Getoor Indices}
As we have pointed out before,  for homogeneous processes like multivariate $\alpha$-stable processes \cite{blumenthalget60}, more general Lévy processes \cite{blumenthalget61}, Feller processes satisfying some mild conditions \cite{schilling98} and homogeneous diffusion with jumps \cite{generalizedindices}, there exists a set of indices, called Blumenthal-Getoor indices, which utilize the symbol to derive maximal inequalities. For a historical overview we refer to \cite{hope}.  Equally, we now want to consider maximal inequalities for non-homogeneous processes with the help of the time-dependent symbol. In this framework, the non-homogeneous growth and sector condition of the symbol play an important role:
\begin{equation}
\lVert q(s, x, \xi) \rVert \leq c(1+\lVert \xi \rVert^2), \tag{IG}\label{IG}
\end{equation}
\begin{equation}
| \Im(q(s,x,\xi) | \leq c_0 \  \Re(p(s,x,\xi)) \tag{IS}\label{IS}
\end{equation}
for every $s\geq 0, x,\xi \in \bbr^d$ and $c,c_0 >0$. \\
Specifically, we want to use the maximal inequalities to examine the paths of the process, including the asymptotic behavior of the sample paths,  the p-variation of the paths and the existence of the exponential moments of the process. The following indices generalize the Blumenthal-Getoor indices as defined in  Definitions 4.2 and 4.5 of \cite{schilling98} or Definition 3.8 of \cite{generalizedindices} to the non-homogeneous case. 

\begin{definition}
The quantities 
\begin{align*}
\beta_0&:=\sup \left\{\lambda \geq 0 : \limsup_{R\to\infty} R^\lambda H(R) =0 \right \} \\
\underline{\beta_0}&:=\sup \left\{\lambda \geq 0 : \liminf_{R\to\infty} R^\lambda H(R) =0 \right \} \\
\overline{\delta_0}&:=\sup \left\{\lambda \geq 0 : \limsup_{R\to\infty} R^\lambda h(R) =0 \right \} \\
\delta_0&:=\sup \left\{\lambda \geq 0 : \liminf_{R\to\infty} R^\lambda h(R) =0 \right \}
\end{align*}
are called \emph{time-dependent indices of $X$ in the starting point},
where 
\begin{align}
H(R)&:= \sup_{s \geq 0} \sup_{y \in \bbr^d} \sup_{\norm{\varepsilon}\leq 1} \abs{p\left(s,y,\frac{\varepsilon}{R}\right)}, \text{ and}\\
h(R)&:= \inf_{s \geq 0} \inf_{y \in \bbr^d} \;  \sup_{\norm{\varepsilon}\leq 1} 
\Re \left( p\left(s,y,\frac{\varepsilon}{4\kappa R} \right)\right)
\end{align}
with $\kappa=(4 \arctan (1/2 c_0))^{-1}$ where $c_0$ comes from the sector condition \eqref{IS}.  
\end{definition}

\begin{definition}
Let $\tau \in \bbr_+, x\in \bbr^d$ and $R>0$. The quantities 
\begin{align*}
\beta_\infty^{\tau,x}&:=\inf \left\{\lambda > 0 : \limsup_{R\to 0} R^\lambda H(\tau,x,R) =0 \right \} \\
\underline{\beta_\infty^{\tau,x}}&:=\inf \left\{\lambda > 0 : \liminf_{R\to 0} R^\lambda H(\tau,x,R) =0 \right \} \\
\overline{\delta_\infty^{\tau,x}}&:=\inf \left\{\lambda > 0 : \limsup_{R\to 0} R^\lambda h(\tau,x,R) =0 \right \} \\
\delta_\infty^{\tau,x}&:=\inf \left\{\lambda > 0 : \liminf_{R\to 0} R^\lambda h(\tau,x,R) =0 \right \}
\end{align*}
are the \emph{time-dependent indices of $X$ at infinity}, where 
\begin{align}
H(\tau,x,R)&:= \sup_{\norm{y-x}\leq 2R} \sup_{\norm{\varepsilon}\leq 1} \abs{p\left(\tau,y,\frac{\varepsilon}{R}\right)}, \text{ and}\\
h(\tau, x,R)&:= \inf_{\norm{y-x}\leq 2R} \;  \sup_{\norm{\varepsilon}\leq 1} 
\Re \left( p\left(\tau,y,\frac{\varepsilon}{4\kappa R} \right)\right)
\end{align}
with $\kappa=(4 \arctan (1/2 c_0))^{-1}$ where $c_0$ comes from the sector condition \eqref{IS}.  
\end{definition}

The proofs of the previous section crucially rely on the  space-time process to transfer properties of homogeneous Markov processes to the non-homogeneous framework.
However, when deriving properties like the asymptotic behavior of sample paths or maximal inequalities we do not expect the space-time process to be of much use. This is due to the fact that if we utilize the space-time process,  a deterministic drift with slope 1 is added. This obscures the path-behavior of the original process. 

In this section we assume all characteristics to be encountered to be with respect to the truncation function $h=id\cdot \chi$, where $\chi \in C_c^\infty(\bbr^d)$ is a symmetric cut-off function with 
$$
1_{B_r(0)} \leq \chi \leq 1_{B_{2R}(0)} 
$$
where $R >0$.

\begin{theorem} \label{prop:techmain}
Let $X$ be a non-homogeneous Itô process with characteristics as in Theorem \ref{main:symbol}. In this case we have
\begin{align} \label{firstestimate}
\bbp^{\tau,x} \left(\sup_{\tau \leq s \leq \tau+t} \lVert X_s- x \rVert \geq R \right) \leq c_d \cdot t \cdot \sup_{\tau < s \leq \tau+t} H(s,x,R)
\end{align}
for $t\geq 0$, $R>0$ and a constant $c_d>0$ which only depends on the dimension $d$. 
If, in addition, \eqref{IS} holds true, we have 
\begin{align} \label{secondestimate}
\bbp^{\tau,x} \left(\sup_{\tau \leq s \leq \tau+t} \lVert X_s- x \rVert < R \right) \leq c_k \cdot \frac{1}{t} \cdot  \frac{1}{\inf_{\tau < s \leq \tau+t} \ h(s,x,R)}
\end{align}
where $c_k>0$ only depends on $c_0$ of the sector condition \eqref{IS}.
\end{theorem}

\begin{remark}
Throughout the following proof we often make use of Lemma 5.2 of \cite{generalizedindices} although it is a statement for time-homogeneous symbols.  A closer look shows that an analogous statement holds for the time-dependent probabilistic symbol and the proof also works analogously.
\end{remark}

\begin{proof}
The proof of \eqref{firstestimate} closely follows the proof of Proposition 3.10 of \cite{generalizedindices}. We omit the proof of \eqref{secondestimate} since it is generalized from the proof of Lemma 6.3 of \cite{schilling98} analogously to the following generalization.  However, let us mention that one has to use Dynkin's formula for non-homogeneous processes as stated in \cite{rueschendorf} when Corollary 3.6 is utilized in \cite{schilling98}.  \\ 
In order to prove \eqref{firstestimate} let $X$ be a non-homogeneous Itô process such that the differential characteristics $(\ell,Q,\nu)$ of $X$ are locally bounded and continuous. At first we show that for $S,R$ and $\sigma:= \inf \{ t \geq \tau: \Vert X_t -x \Vert > S \}$ as above we have
\begin{align} \label{techmainstopped}
\bbp^{\tau,x} \left(\sup_{\tau \leq s \leq t} \lVert X_s^\sigma-x \rVert  \geq 2R \right) \leq c_d \cdot (t-\tau) \cdot\sup_{\tau \leq s \leq t}  \sup_{\norm{y-x}\leq S} \sup_{\norm{\varepsilon}\leq 1} \abs{p\left(s,y,\frac{\varepsilon}{2R} \right) }
\end{align}
where $c_d=4d+16\widetilde{c}_d$ and $t \geq \tau$.
We introduce the stopping time 
\[
\tau_R:= \inf \{t\geq \tau: \norm{\Delta X_t^\sigma} >R \},
\]
as the first time the jumps of $X^\sigma$ exceed $R$, and estimate the following
\begin{align} \label{separateone}
\bbp^{\tau,x} \left(\sup_{\tau \leq s \leq t} \lVert X_s^\sigma-x \rVert  \geq 2R \right) \leq \bbp^{\tau,x} \left(\sup_{\tau \leq s \leq t} \lVert X_s^\sigma-x \rVert  \geq 2R, \tau_R > t \right) + \bbp^{\tau,x} \Big(\tau_R \leq t \Big).
\end{align}
We deal with the terms on the right-hand side one after another, starting with the first one.

Again we separate the first term \eqref{separateone} in order to get control over the big jumps.  For $t \geq \tau$ let 
$$
\check{X}_t:= X_t- \sum_{\tau \leq s \leq t} \Delta X_s(1-\chi(\Delta X_s)).
$$
The process $\check{X}$ is a special semimartingale on $[\tau, \infty)$ with characteristics 
\begin{align*}
\check{B}_t^{(i)}&= \int_\tau^t b^i(s,X_s) \ ds, \quad \bbp^{\tau,x}\text{-a.s.} \\
\check{C}_t^{(ij)} &= \int_\tau^t c^{ij}(s,X_s) \ ds, \quad \bbp^{\tau,x}\text{-a.s.}\\
\check{\nu}( ; dt,dx) &= \chi(y) 1_{\csi{0}{\sigma}} (t) \ N(X_t,dy) \ dt, \quad \bbp^{\tau,x}\text{-a.s.}
\end{align*}
Now let $u=(u_1,...,u_d)':\bbr^d \to \bbr^d$ be such that $u_j \in C_b^2(\bbr^d)$ is 1-Lipschitz continuous, $u_j$ depends only on $x^{(j)}$ and is zero in zero for $j=1,...,d$. 
We define the auxiliary process for $t \geq \tau$:
\[
\check{M}_t:=u(\check{X}_t^\sigma-x) - \int_\tau^{t\wedge \sigma} \sum_{j=1}^d F^{(j)}_s \ ds
\]
where 
\begin{align} \begin{aligned} \label{auxprocess}
F_s^{(j)} &= \partial_j u(\check{X}_{s-}-x) \ell^{(j)}(s,X_{s-}) -\frac{1}{2} \partial_j \partial_j u(\check{X}_{s-}-x) Q^{(jj)} (s,X_{s-}) \\
  &-\int_{z\neq 0} \Big( u(\check{X}_{s-}-x+z)-u(\check{X}_{s-}-x) - \chi(z) z^{(j)} \partial_j u(\check{X}_{s-}-x)\Big) \chi(z) \ N_s(X_{s-}, dz).
\end{aligned} \end{align}
$\check{M}$ is a $\bbp^{\tau,x}$-local martingale on $[\tau,\infty)$ by \cite{jacodshir} Theorem II.2.42. Applying Lemma 3.7 of \cite{mydiss} we have under \eqref{IG} 
\[
\abs{F_s^{(j)}} \leq const \cdot \sum_{0\leq \abs{\alpha}\leq 2} \norm{\partial^\alpha u}_\infty
\]
since $u_j\in C_b^2(\bbr^d)$.  Let us mention, that although Lemma 3.7 of \cite{mydiss} considers homogeneous diffusion with jumps only the proof is alike for non-homogeneous Itô processes. In particular,  since $\check{M}$ is uniformly bounded it is a $L^2$-martingale on $[\tau,t]$.
We define 
\[
  D:=\left\{\omega\in\Omega:\int_\tau^{t\wedge \sigma(\omega)} \norm{F_s(\omega)} \ ds \leq R \right\}
\]
and obtain for $t \geq \tau$:
\begin{align}\label{separatetwo}
\bbp^{\tau,x}  \left(\sup_{\tau \leq s \leq t} \lVert X_s^\sigma-x \rVert \geq 2R, \tau_R > t \right) \leq \bbp^{\tau,x} \left(\sup_{\tau \leq s \leq t} \lVert X_s^\sigma-x \rVert \geq 2R, \tau_R > t, D \right) + \bbp^{\tau,x} (D^c).
\end{align}
For the left summand of \eqref{separatetwo} we estimate for again for $t \geq \tau$:
\begin{align*}
\bbp^{\tau,x} \left( \sup_{\tau \leq s \leq t} \lVert u(X_s^\sigma- x) \rVert  \geq 2R, \tau_R > t, D \right) 
&= \bbp^{\tau,x} \left( \sup_{\tau \leq s \leq t} \lVert u(\check{X}_s^\sigma- x) \rVert  \geq 2R, \tau_R > t, D \right) \\
&\leq \bbp^{\tau,x} \left( \sup_{\tau \leq s \leq t} \lVert u(X_s^\sigma- x) \rVert -\int_\tau^{t \wedge \sigma} F_s \ ds \geq R, \tau_R > t, D \right) \\
&\leq \bbp^{\tau,x}\left(\sup_{\tau \leq s \leq t} \lVert \check{M}_{s \wedge \sigma} \rVert \geq R \right) \\
&= \bbp^{\tau,x} \left(\sup_{0 \leq s \leq t-\tau} \lVert \check{M}_{(s+\tau) \wedge \sigma}  \rVert \geq R \right) \\
&\leq \frac{1}{R^2} \bbe^{\tau,x} \left(\norm{\check{M}_t^\sigma}^2  \right) \\
&\leq \frac{1}{R^2} \sum_{j=1}^d  \bbe^{\tau,x} \left( \left[\check{M}^{(j)}, \check{M}^{(j)} \right]^\sigma_t \right) \\
&\leq \frac{1}{R^2} \sum_{j=1}^d \bbe^{\tau,x} \left( [\check{X}_\cdot^{(j)}, \check{X}_\cdot^{(j)}]_t^\sigma   \right),
\end{align*} 
where we used Doob's inequality for the martingale $\check{M}^\sigma$ and the Lipschitz property of $u$ in combination with Corollary II.3 of \cite{protter}.
Since 
\begin{align*}
\bbe^{\tau,x} \left( [\check{X}_\cdot^{(j)}, \check{X}_\cdot^{(j)}]_t^\sigma  \right) 
= \bbe^{\tau,x} \left( \left<\check{X}_\cdot^{(j),c}, \check{X}_\cdot^{(j),c}\right>_t^\sigma  \right)
+ \bbe^{\tau,x} \left( \int_\tau^{t\wedge \sigma} \int_{z\neq 0} (z^{(j)})^2 \chi(z)^2 \ N_s(X_s,dz) \ ds\right).
\end{align*}
we obtain 
\begin{align*}
&\bbp^{\tau,x} \Big( \sup_{\tau \leq s \leq t} \lVert u(X_s^\sigma- x) \rVert \geq 2R, \tau_R > t, D \Big) \\
& \leq \frac{1}{R^2} \sum_{j=1}^d \bbe^{\tau,x} \int_\tau^{t\wedge \sigma} Q^{(jj)}(s,X_s) \ ds 
  +  \bbe^{\tau,x} \int_\tau^{t\wedge \sigma} \int_{z\neq 0} \frac{\norm{z}^2}{R^2} \chi(z)^2 \ N_s(X_s,z) \ ds \\
& \leq 4\cdot  \sum_{j=1}^d \bbe^{\tau,x} \int_\tau^{t\wedge \sigma} \left(\frac{e_j}{2R} ' Q(s,X_s) \frac{e_j}{2R}\right) \ ds + 4^2  \bbe^{\tau,x} \int_\tau^{t\wedge \sigma} \int_{z\neq 0} \left( \norm{\frac{z}{2R}}^2 \wedge 1 \right) \ N_s(X_s,dz) \ ds \\
& \leq 4\cdot  \sum_{j=1}^d \bbe^{\tau,x} \sup_{\tau< s<t\wedge \sigma} \left(\frac{e_j}{2R} ' Q(s,X_s) \frac{e_j}{2R}\right) \ ds + 4^2  \bbe^{\tau,x} \sup_{\tau< s<t\wedge \sigma} \int_{z\neq 0} \left( \norm{\frac{z}{2R}}^2 \wedge 1 \right) \ N_s(X_s,dz)  \\
& \leq 4(t-\tau)\cdot  \sum_{j=1}^d \sup_{\tau< s<t\wedge \sigma} \sup_{\norm{y-x} \leq S} \Re p\left(s,y, \frac{e_j}{2R} \right) + 
   4^2(t-\tau) \sup_{\tau< s<t\wedge \sigma} \sup_{\norm{y-x}\leq S} \int_{z\neq 0} \left( \norm{\frac{z}{2R}}^2 \wedge 1 \right) \ N_s(y,dz)  \\ 
& \leq 4(t-\tau)d \sup_{\tau < s<t\wedge \sigma} \sup_{\norm{y-x} \leq S} \sup_{\norm{\varepsilon} \leq 1} \abs{p\left(s,y,\frac{\varepsilon}{2R} \right) } + 
   4^2 t \sup_{\tau <s<t\wedge \sigma} \sup_{\norm{y-x}\leq S} \widetilde{c}_d \sup_{\norm{\varepsilon} \leq 1}
   \abs{p\left(s,y,\frac{\varepsilon}{2R}\right)} 
\end{align*} \normalsize
where we have used Lemma 5.2 of \cite{generalizedindices} on the second term.  By choosing a seqence $(u_n)_{n\in\bbn}$ of functions of the type described above which tends to the identity in a monotonous way we obtain
\begin{align} \label{finalterm1D}
\bbp^{\tau,x} \left( \sup_{\tau \leq s \leq t} \lVert X_s^\sigma- x \rVert  \geq 2R, \tau_R > t, D \right) \leq 
(4d+ 4^2 \widetilde{c}_d ) (t-\tau) \sup_{\tau <s<t\wedge \sigma}\sup_{\norm{y-x} \leq S} \sup_{\norm{\varepsilon} \leq 1} \abs{p\left(s,y,\frac{\varepsilon}{2R} \right) } .
\end{align}
Now let us consider the term $\bbp^{\tau,x}$ of \eqref{separatetwo}. The Markov inequality provides
\begin{align*}
\bbp^{\tau,x}(D^c)=\bbp^{\tau,x} \left( \int_\tau^{t\wedge \sigma} \norm{F_s} \ ds > R \right) 
\leq \frac{1}{R} \sum_{j=1}^d \bbe^{\tau,x} \Big(\int_\tau^{t\wedge \sigma} \abs{F_s^{(j)}}  \ ds\Big)=:(*)
\end{align*} 
Again we chose a sequence $(u_n)_{n\in\bbn}$ of functions as we described in \eqref{auxprocess}, but this time it is important that the first and second derivatives are uniformly bounded. Since the $u_n$ converge to the identity, the first partial derivatives tend to 1 and the second partial derivatives to 0. In the limit ($n\to\infty$) we obtain
\begin{align}
(*)&\leq \frac{1}{R} \sum_{j=1}^d \bbe^{\tau,x} \int_\tau^{t \wedge \sigma} \abs{\ell^{(j)}(s,X_s) +
   \int_{z\neq 0} (-z^{(j)} \chi(z) + (\chi(z))^2 z^{(j)} ) \ N_s(X_s,dz) } \ ds \nonumber  \\
   &\leq 2\sum_{j=1}^d \bbe^{\tau,x} \int_\tau^{t\wedge \sigma} \abs{\frac{\ell^{(j)}(s,X_s)}{2R} + 
   \int_{z\neq 0} \sin\left(\frac{z'e_j}{2R} \right) - \frac{z^{(j)}\chi(z)}{2R} \ N_s(X_s,dz)} \ ds  \label{termA}\\
   &+2\sum_{j=1}^d \bbe^{\tau,x} \int_\tau^{t\wedge\sigma} \abs{\int_{z\neq 0} \frac{(\chi(z))^2 z^{(j)}}{2R} - \sin\left(\frac{z'e_j}{2R} \right) \ N(X_s,dz) } \ ds. \label{termD}
\end{align}
For term \eqref{termA} we get
\begin{align}
&2 \sum_{j=1}^d \bbe^{\tau,x} \int_\tau^{t\wedge \sigma} \abs{\frac{\ell(s,X_s)'e_j}{2R} 
  + \int_{z\neq 0} \sin\left(\frac{z'e_j}{2R} \right) - \frac{z'e_j \chi(z)}{2R} \ N_s(X_s,dz)} \ ds \nonumber \\
&\leq 2(t-\tau)\sum_{j=1}^d \sup_{\tau < s\leq t\wedge \sigma} \bbe^{\tau,x} \abs{\frac{\ell(s,X_s)'e_j}{2R} 
  + \int_{z\neq 0} \sin\left(\frac{z'e_j}{2R} \right) - \frac{z'e_j \chi(z)}{2R} \ N_s(X_s,dz)} \nonumber \\
&\leq 2(t-\tau)d \sup_{\tau < s\leq t}\sup_{\norm{y-x}\leq S} \sup_{\norm{\varepsilon}\leq 1} \abs{\Im p\left(s,y,\frac{\varepsilon}{2R} \right)} \label{finalterm1DcA}
\end{align}
and for term \eqref{termD}
\begin{align}
&2\sum_{j=1}^d \bbe^{\tau,x} \int_\tau^{t\wedge \sigma} \abs{ \int_{z\neq 0} \frac{(\chi(z))^2 z'e_j}{2R}-
\sin\left(\frac{z'e_j}{2R} \right) \ N_s(X_s,dz)} \ ds \nonumber \\
&\leq 2\sum_{j=1}^d \bbe^{\tau,x} \int_\tau^{t\wedge \sigma} \abs{\int_{B_{2R}(0)\backslash\{0\}} 1- \cos\left(\frac{z'e_j}{2R} \right) \ N_s(X_s,dz) }\nonumber \\
& \hspace{10mm}+\abs{\int_{B_{2R}(0)^c} 1 \ N_s(X_s,dz)} \ ds \nonumber \\
&\leq 2(t-\tau)d \sup_{\tau < s\leq t} \sup_{\norm{y-x}\leq S} \sup_{\norm{\varepsilon} \leq 1} \Re p\left(s,y,\frac{\varepsilon}{2R} \right) +  2^2 (t-\tau)d \sup_{\tau < s\leq t} \sup_{\norm{y-x}\leq S} \widetilde{c}_d \sup_{\norm{\varepsilon} \leq 1}
   \abs{p\left(s,y,\frac{\varepsilon}{2R}\right)} \label{finalterm1DcB}
\end{align}
where we have used again Lemma 5.2 of \cite{generalizedindices} on the second term.

It remains to deal with the second term of \eqref{separateone}.
Let $\delta >0$ be fixed (at first) and $m:\bbr \to ]1,1+\delta[$ a strictly monotone increasing auxiliary function. Since $m\geq 1$ and since we have at least one jump of size $>R$ on $\{\tau_R\leq t\}$ we obtain
\begin{align*}
\bbp^{\tau,x}(\tau_R\leq t) & \leq \bbp^{\tau,x} \left( \int_\tau^t \int_{\norm{z}\geq R} m(\norm{z}) \ \mu^{X^\sigma}(\cdot;ds,dz) \geq m(R) \right) \\
& \leq \frac{1}{m(R)} \bbe^{\tau,x} \left( \int_0^t \int_{\norm{z}\geq R} m(\norm{z}) 1_{\csi{0}{\sigma}}(s) \ \mu^X(\cdot;ds,dz) \right) \\
& = \frac{1}{m(R)} \bbe^{\tau,x} \left( \int_\tau^t \int_{z\neq 0} m(\norm{z}) 1_{\rosi{0}{\sigma}}(s) 1_{B_R(0)^c} (z) \ N_s(X_s,dz) \ ds \right) \\
&\leq (1+\delta) (t-\tau) \sup_{\tau < s\leq t\wedge \sigma} \bbe^{\tau,x} (N_s(X_s, B_R(0)^c)) \\
&\leq (1+\delta) (t-\tau) \sup_{\tau < s\leq t} \sup_{\norm{y-x}\leq S} N_s(y,B_R(0)^c) \\
&\leq (1+\delta) 4(t-\tau) \sup_{\tau < s\leq t} \sup_{\norm{y-x}\leq S} \int_{z\neq 0} \left(\norm{\frac{z}{2R}}^2 \wedge 1 \right) \ N_s(y,dz)
\end{align*}
because $m(\norm{z}) 1_{\rosi{0}{\sigma}}(s) 1_{B_R(0)^c} (z)$ is in class $F_p^1$  of Ikeda and Watanabe (see \cite{ikedawat}, Section II.3). Since $\delta$ can be chosen arbitrarily small we obtain by Lemma 5.2 of \cite{generalizedindices}
\begin{align} \label{finalterm2}
\bbp^{\tau,x}(\tau_R\leq t)  \leq  4(t-\tau) \sup_{\tau < s\leq t} \sup_{\norm{y-x}\leq S} \widetilde{c}_d \sup_{\norm{\varepsilon} \leq 1} \abs{p\left(s,y,\frac{\varepsilon}{2R}\right)}.
\end{align}

Plugging together \eqref{finalterm1D}, \eqref{finalterm1DcA}, \eqref{finalterm1DcB} and \eqref{finalterm2} we obtain \eqref{techmainstopped}.  For the particular case $\sigma=\sigma_{3\widetilde{R}}^x$ we have
\[
\left\{ \sup_{\tau \leq s \leq t} \lVert u(X_s^\sigma- x) \rVert \geq 2\widetilde{R} \right\} = \left\{ \sup_{\tau \leq s \leq t} \lVert u(X_s- x) \rVert \geq 2\widetilde{R} \right\},
\]
and therefore, for every $\widetilde{R}>0$
\begin{align} \label{techmain}
\bbp^{\tau,x}\left(\sup_{\tau \leq s \leq t} \lVert u(X_s^\sigma- x) \rVert \geq 2\widetilde{R}\right) \leq c_d \cdot (t-\tau) \cdot \sup_{\tau < s\leq t} \sup_{\norm{y-x}\leq 3\widetilde{R}} \sup_{\norm{\varepsilon}\leq 1} \abs{p\left(s,y,\frac{\varepsilon}{2\widetilde{R}} \right) }.
\end{align}
Setting $R:=(1/2)\widetilde{R}$ we obtain \eqref{firstestimate}. 
\end{proof}

The maximal inequalities \eqref{firstestimate} and \eqref{secondestimate} allow for  statements of the asymptotic behavior of the sample paths: First, we state a result concerning the behavior near the starting point $x \in \bbr^d$ at time $\tau \geq 0$ with respect to the measure $\bbp^{\tau,x}$. The second statement, treats the same behavior at infinity.  The proof of both statements is inspired by Theorem 4.3 and Theorem 4.6 of \cite{schilling98} but takes the time-component $\tau$ into account. 
Let us mention, that for a function $f:\bbr_+ \to \bbr_+$ we denote by $f(t+):= \lim_{s \downarrow t} f(s)$, if the limit exists.

\begin{theorem} \label{thm:tto0}
Let $X$ be a non-homogeneous Itô process such that the differential characteristics of $X$ are locally bounded and continuous. Then we have
\begin{align}
\lim_{t\to 0} \ t^{-1/\lambda}\sup_{\tau \leq s \leq \tau+t} \lVert X_s-x \rVert &=0  \text{ for all } \lambda > \beta_\infty^{\tau+,x} \label{proved1} \\
\liminf_{t\to 0 } \ t^{-1/\lambda}\sup_{\tau \leq s \leq \tau+t} \lVert X_s-x \rVert &=0  \text{ for all } \beta_\infty^{\tau+, x} \geq \lambda > \underline{\beta_\infty^{\tau+,x}}.\label{unproved2}
\end{align}
If the symbol $p(\tau,x,\xi)$ of the process $X$ satisfies \eqref{IS} then we have in addition
\begin{align}
\limsup_{t\to 0 } \  t^{-1/\lambda}\sup_{\tau \leq s \leq \tau+t} \lVert X_s-x \rVert=\infty & \text{ for all } \overline{\delta_\infty^{\tau +,x}} > \lambda \geq \delta_\infty^{\tau +,x}  \label{proved2} \\
\lim_{t\to 0 }  \  t^{-1/\lambda}\sup_{\tau \leq s \leq \tau+t} \lVert X_s-x \rVert=\infty & \text{ for all } \delta_\infty^{\tau +,x} > \lambda.\label{unproved1}
\end{align}
All these limits are meant $\bbp^{\tau,x}$-a.s with respect to every $\tau \in \bbr_+$ and $x\in\bbr^d$.
\end{theorem}

\begin{proof}
Here, we only prove \eqref{proved1} and \eqref{proved2} and omit the proofs of \eqref{unproved2} and \eqref{unproved1} since \eqref{proved1} and \eqref{unproved1} and \eqref{unproved2} and \eqref{proved2} are very similar, respectively.\\ 
Let $\varepsilon >0$, $\tau \in \bbr_+$ and $x\in\bbr^d$.  We start with proving \eqref{proved1}:\\
Let $\lambda> \sup_{\tau < s \leq \tau + \varepsilon} \beta_\infty^{s,x}$ and choose $\sup_{\tau < s \leq \tau + \varepsilon} \beta_\infty^{s,x} < \alpha_2 < \alpha_1 <\lambda$.  For $t<T_0^\varepsilon$ with $T_0^\varepsilon$ sufficiently small \eqref{firstestimate} one obtains:
\begin{align*}
\bbp^{\tau,x}\left(\sup_{\tau \leq s \leq \tau+t} \lVert X_s-x \rVert \geq t^{1/\alpha_1}\right) \leq c_d \cdot t \cdot \sup_{\tau < s\leq \tau+\varepsilon} H(s,x,t^{1/\alpha_1}) 
\leq c_d' \cdot t (t^{1/\alpha_1})^{-\alpha_2}
=c_d' t^{1-(\alpha_2/\alpha_1)}.
\end{align*}
Now let $t_k:=(1/2)^k$ for $k\in\bbn$. Since 
\[
\sum_{k=k_0 ^\varepsilon}^\infty \bbp^{\tau,x} \left( \sup_{\tau \leq s \leq \tau+t_k} \lVert X_s-x \rVert \geq (t_k)^{1/\alpha_1} \right) \leq c_d' \sum_{k=k_0^\varepsilon}^\infty 2^{-k(1-(\alpha_2/\alpha_1))} <+\infty
\]
where $k_0^\varepsilon$ depends on $T_0^\varepsilon$, the Borel-Cantelli Lemma derives
\[
\bbp^{\tau,x} \left( \limsup_{k\to\infty} \sup_{\tau \leq s \leq \tau+t_k} \lVert X_s-x \rVert \geq (t_k)^{1/\alpha_1}\right) =0,
\]
and, hence, $\sup_{\tau \leq s \leq \tau+t_k} \lVert X_s-x \rVert < (t_k)^{1/\alpha_1}$ for all $k\geq k_1^\varepsilon(\omega)$ on a set of probability one. For fixed $\omega$ in this set and $t \in [t_{k+1}, t_k]$ and $k\geq k^\varepsilon_1(\omega) \geq k^\varepsilon_0$, we have 
\[
\sup_{\tau \leq s \leq \tau+t} \lVert X_s(\omega)-x \rVert  \leq \sup_{\tau \leq s \leq \tau+t_k} \lVert X_s(\omega)-x \rVert \leq t_k^{1/\alpha_1} \leq 2^{1/\alpha_1} t^{1/\alpha_1}
\]
and since $\lambda >\alpha_1$
\[
t^{-1/\lambda} \sup_{\tau \leq s \leq \tau+t} \lVert X_s-x \rVert  \leq 2^{1/\alpha_1} t^{(1/\alpha_1)-(1/\lambda)}
\]
which converges $\bbp^{\tau,x}$-a.s to zero as $t\downarrow 0$.  Since $\lambda>\sup_{\tau < s \leq \tau + \varepsilon} \beta_\infty^{s,x}$ and $\varepsilon>0$ arbitrary, $\varepsilon \downarrow 0$ provides the statement.\\ 
In order to prove \eqref{proved2} we derive the following: \\
Let $ \inf_{\tau < s \leq \tau+\varepsilon} \overline{\delta_\infty^{s,x}}> \lambda' > \lambda$.  Moreover, let $(t_k)_{k \in \bbn}$ be a sequence such that 
$$
\lim_{k \to \infty} (t_k)^{\lambda'} h(s,x,t_k)= \infty, \quad \forall s \in ]\tau, \tau+\varepsilon].
$$
Hence, the maximal inequality \eqref{secondestimate} provides for $k$ sufficiently large
\begin{align*}
 \bbp^{\tau,x} \left(\sup_{\tau \leq s \leq \tau+(t_k)^{\lambda'}} \lVert X_s- x \rVert < t_k \right)& \leq c_d \cdot \frac{1}{(t_k)^{\lambda'}} \cdot \frac{1}{\inf_{\tau < s \leq \tau+(t_k)^{\lambda'}} h(s,x,t_k)} \\
&\leq c_d \cdot \frac{1}{(t_k)^{\lambda'}} \cdot \frac{1}{\inf_{\tau < s \leq \tau+\varepsilon} h(s,x,t_k)} \\
&\underset{k \to \infty}\longrightarrow 0.
\end{align*}
Fatou's lemma implies 
\begin{align*}
0 &= \liminf_{k \to \infty} \bbp^{\tau,x} \left(\sup_{\tau \leq s \leq \tau+(t_k)^{\lambda'}} \lVert X_s- x \rVert < t_k \right)\\
&= 1- \limsup_{k \to \infty} \bbp^{\tau,x} \left(\sup_{\tau \leq s \leq \tau+(t_k)^{\lambda'}} \lVert X_s- x \rVert \geq t_k \right)\\
&\geq 1-  \bbp^{\tau,x} \left(\limsup_{k \to \infty} \left\{\sup_{\tau \leq s \leq \tau+(t_k)^{\lambda'}} \lVert X_s- x \rVert \; \geq t_k \right\} \right).
\end{align*}
Hence, 
$$
\bbp^{\tau,x} \left( \left(\frac{1}{t_k}\right)\sup_{\tau \leq s \leq \tau+(t_k)^{\lambda'}} \lVert X_s- x \rVert \geq 1 \text{, infinitely often}\right)=1,
$$
and, therefore,
$$
\limsup_{k \to \infty} \left( \frac{1}{t_k}\right)^{-\frac{1}{\lambda'}}\sup_{\tau \leq s \leq \tau+t_k} \lVert X_s- x \rVert \geq 1.
$$
Since $\lambda < \lambda'$ we observe, that
$$
\limsup_{t \to 0} \left( \frac{1}{t}\right)^{-\frac{1}{\lambda}}\sup_{\tau \leq s \leq \tau+t} \lVert X_s- x \rVert = \infty.
 $$
Since $\varepsilon >0$ arbitrary, the statement follows.
\end{proof}

The proof of the following lemma parallels Lemma \ref{thm:tto0}, and, hence, we omit details of the proof.

\begin{theorem} \label{thm:ttoinfty}
Let $X$ be a non-homogeneous Itô process such that the differential characteristics of $X$ are locally bounded and continuous. Then we have
\begin{align}
\lim_{t\to \infty} \ t^{-1/\lambda}\sup_{\tau \leq s \leq \tau+t} \lVert X_s-x \rVert &=0  \text{ for all } \lambda < \beta_0 \\
\liminf_{t\to \infty} \ t^{-1/\lambda}\sup_{\tau \leq s \leq \tau+t} \lVert X_s-x \rVert &=0  \text{ for all } \underline{\beta_0} > \lambda \geq \beta_0
\end{align}
If the symbol $p(\tau,x,\xi)$ of the process $X$ satisfies \eqref{IS} then we have in addition
\begin{align}
\limsup_{t\to \infty} \  t^{-1/\lambda}\sup_{\tau \leq s \leq \tau+t} \lVert X_s-x \rVert=\infty & \text{ for all } \overline{\delta_0} < \lambda \leq \delta_0 \\
\lim_{t\to \infty}  \  t^{-1/\lambda}\sup_{\tau \leq s \leq \tau+t} \lVert X_s-x \rVert=\infty & \text{ for all } \delta_0 < \lambda.
\end{align}
All these limits are meant $\bbp^{\tau,x}$-a.s with respect to every $\tau \in \bbr_+$ and $x\in\bbr^d$.
\end{theorem}

The time homogeneous version of the following result is Lemma 5.11 of \cite{levymatters3}.

\begin{lemma}
Let $X$ be a non-homogeneous Itô process such that the differential characteristics of $X$ are locally bounded and continuous.  Moreover, we have 
\begin{equation}\label{boundednesscond}
\sup_{s\geq 0, x \in \bbr^d} |p(s,x,\xi) | \leq c(1+\lVert \xi \rVert^2) \quad \text{for }\xi \in \bbr^d.
\end{equation}
Then 
$$\bbe^{\tau,x}\left( e^{X_t'\xi} \right)< \infty$$
 for all $t \geq \tau\geq 0$ and $x,\xi \in \bbr^d$.
\end{lemma}

\begin{proof}
At first, let us reconsider the stopping time $\sigma$ as defined in \eqref{exittimeR}: 
$$
\sigma_R^{\tau,x} := \sigma := \inf\{ h \geq \tau: \Vert X^{\tau,x}_{h} -x \Vert > R\}
$$
and let $p(\tau,x,\xi)$ be the symbol of $X$ given by 
$$
p(\tau,x,\xi) = -i \ell(\tau,x)'\xi + \frac{1}{2} \xi' Q(\tau,x) \xi - \int_{y \neq 0} \left(e^{iy' \xi} -1 - iy' \xi \cdot \chi(y) \right) \ N_\tau(x,dy).
$$
In order to apply Gronwall's Lemma we estimate the following: 
\begin{align*}
\bbe^{\tau,x} \left( e^{(X_t^{\sigma}-x)\xi} \right)&= \bbe^{\tau,x} \left[ \int_{\tau}^{t\wedge \sigma} e^{(X^\sigma_{s}-x)\xi} \left(\xi \ell(s,X_s)-  \frac{1}{2}\xi^2\ Q(s,X_s) \right.\right.\\
&\quad + \left. \left.\int_{ \{y\neq 0\}}  (e^{\xi y}-1 - \xi y\chi(y)) \ N_s(X_s,dy) \right)\ ds\right]\\
&\leq b(\xi) \int_\tau^t \bbe^{\tau,x}  \left( e^{(X_s^{\sigma}-x)\xi} \right) \ ds,
\end{align*}
Condition \eqref{boundednesscond} provides (see Lemma 3.3 of \cite{mydiss}) the finiteness of the constant 
$$
b(\xi):=\sup_{s\geq 0, x \in \bbr^d} \left|\xi \ell(s,x)-  \frac{1}{2}\xi^2\ Q(s,x) + \int_{ \{y\neq 0\}}  (e^{\xi y}-1 - \xi y\chi(y)) \ N_s(x,dy) \right|. 
$$
 Now, application of Gronwall's Lemma provides 
\begin{align*}
\bbe^{\tau,x}  \left( e^{(X_s^{\sigma}-x)\xi} \right) &\leq 1+ b(\xi) \int_\tau^t e^{b(\xi)(t-s)} \ ds = e^{b(\xi)(t-\tau)},
\end{align*}
and with Fatou's Lemma follows
$$
\bbe^{\tau,x}  \left( e^{(X_s-x)\xi} \right) \leq \liminf_{R \to \infty} \bbe^{\tau,x}  \left( e^{(X_s^{\sigma}-x)\xi} \right) \leq e^{b(\xi)(t-\tau)},
$$
where the first inequality follows by Proposition \ref{prop:techmain}.
\end{proof}

Finally we generalize Theorem 2.10 of \cite{manstaviciusschnurr} which is a nice and applicable criterion for the finiteness of p-variation. 

\begin{theorem}\label{p-Var}
Let $X$ be a non-homogeneous Itô process such that the differential characteristics of $X$ are locally bounded and continuous. For every $t \geq \tau \geq 0$ and $p > \sup_{\tau \geq 0} \beta_\infty^{\tau,x}$ the following holds true 
\begin{equation}
\sup_{\pi_n} \sum_{j=1}^n \lVert X_{t_j}-X_{t_{j-1}} \rVert^p < + \infty, \quad \bbp^{\tau,x}\text{-a.s.}
\end{equation}
where the supremum is taken over all finite partitions $\pi_n= (t_i)_{i=1,...,n}$ with $\tau =t_0 < t_1<...<t_n=t$. I.e.,  the p-variation of the paths of $X$ are $\bbp^{\tau,x}$-a.s finite for $p > \sup_{\tau \geq 0} \ \beta_\infty^{\tau,x}$.
\end{theorem}

\begin{proof}
Let $t, r >0$ and $\lambda>p$. Then Theorem \ref{techmain} provides the following estimation 
\begin{align*}
\alpha(t,r)&:= \sup\left\{\bbp^{\tau,x}(\lVert X_s-x \rVert \geq r) : \tau \geq 0,  s \in [\tau, \tau+t], x \in \bbr^d \right\} \\
&\leq \sup\left\{\bbp^{\tau,x}\left(\sup_{\tau \leq s \leq \tau+t} \lVert X_s-x \rVert \geq r\right) : \tau \geq 0,  x \in \bbr^d \right\}\\
&\leq c_d \cdot t \cdot \sup_{\tau \geq 0} \cdot  \sup_{\tau \leq s \leq \tau+t} H(s,x,r) \\
&\leq c_d \cdot t \cdot \sup_{\tau \geq 0} H(\tau,x,r) \\
&\leq c_d \cdot t \cdot K \cdot r^{-\lambda}
\end{align*}
for $r$ small enough and $K>0$. Hence, Theorem 1.3 of \cite{manstavicius} provides the statement. 
\end{proof}

\bibliographystyle{imsart-number} 
\bibliography{bibliography}       

\end{document}